\def\yes{\if00}
\def\no{\if01}
\def\iftwelvept{\yes}
\def\ifusepsfont{\no}  %Yes if you prefer PostScripy fonts
\def\ifusepdf{\no}     %Yes if you want hypertext links

\iftwelvept
\documentclass[leqno,12pt]{amsart}
\else
\documentclass[leqno]{amsart}
\fi
\usepackage{amssymb}
\usepackage{amscd}
\usepackage{latexsym}
\usepackage{verbatim}
\ifusepdf
\usepackage{hyperref}
\else\fi
\ifusepsfont
\usepackage[T1]{fontenc}
\usepackage{times}
\else\fi

\iftwelvept
\else\fi

\theoremstyle{plain}
\newtheorem{Theorem}{Theorem}[section]
\newtheorem{Proposition}[Theorem]{Proposition}
\newtheorem{Lemma}[Theorem]{Lemma}
\newtheorem{Corollary}[Theorem]{Corollary}
\newtheorem{LangConj}{Lang's conjecture in the absolute form\!\!}
\newtheorem{Claim}{Claim}[Theorem]

\theoremstyle{definition}

\newtheorem{Remark}[Theorem]{Remark}

\renewcommand{\theTheorem}{\arabic{section}.\arabic{Theorem}}
\renewcommand{\theClaim}{\arabic{section}.\arabic{Theorem}.\arabic{Claim}}
\renewcommand{\theequation}{\arabic{section}.\arabic{Theorem}.\arabic{Claim}}

\def\rom{\textup}
\newcommand{\ZZ}{{\mathbb{Z}}}
\newcommand{\QQ}{{\mathbb{Q}}}
\newcommand{\RR}{{\mathbb{R}}}
\newcommand{\CC}{{\mathbb{C}}}

\newcommand{\DD}{{\mathcal{D}}}
\newcommand{\ZDD}{{\overline{\mathcal{D}}}}
\newcommand{\Proof}{{\sl Proof.}\quad}
\newcommand{\trdeg}{\operatorname{tr.deg}}
\newcommand{\Spec}{\operatorname{Spec}}
\newcommand{\Span}{\operatorname{Span}}
\newcommand{\Gal}{\operatorname{Gal}}
\newcommand{\Stab}{\operatorname{Stab}}
\newcommand{\acherncl}{\widehat{{c}}}
\newcommand{\adeg}{\widehat{\operatorname{deg}}}

\newcommand{\QED}{{\unskip\nobreak\hfil\penalty50\quad\null\nobreak\hfil
{$\Box$}\parfillskip0pt\finalhyphendemerits0\par\medskip}}
\newcommand{\rest}[2]{\left.{#1}\right\vert_{{#2}}}
\begin{document}

%%%%%%%%%%%
%% Title %%
%%%%%%%%%%%
\title[A generalization of conjectures of Bogomolov and Lang]%
{A generalization of conjectures \\
of Bogomolov and Lang \\
over finitely generated fields}
\author{Atsushi Moriwaki}
\address{Department of Mathematics, Faculty of Science,
Kyoto University, Kyoto, 606-8502, Japan}
\email{moriwaki@kusm.kyoto-u.ac.jp}
\date{18/August/1999, 4:45PM (JP), (Version 1.0)}

\maketitle

\renewcommand{\theTheorem}{\Alph{Theorem}}

\section*{Introduction}
Let $K$ be a finitely generated field over $\QQ$ with
$d = \trdeg_{\QQ}(K)$, and let
$\overline{B}$ 
be a big polarization of $K$.
Let $A$ be an abelian variety over $K$, and $L$ a symmetric
ample line bundle on $A$.
In the paper \cite{MoArht}, we define
the height pairing
\[
\langle\ , \ \rangle_L^{\overline{B}} :
A(\overline{K}) \times A(\overline{K}) \to \RR
\]
assigned to $\overline{B}$ and $L$
with properties:
$\langle x , x \rangle_L^{\overline{B}} \geq 0$ for all $x \in
A(\overline{K})$ and the equality holds if and only if
$x \in A(\overline{K})_{tor}$.
For $x_1, \ldots, x_l \in A(\overline{K})$,
we denote $\det \left( \langle x_i, x_j \rangle_L^{\overline{B}} \right)$
by $\delta_L^{\overline{B}}(x_1, \ldots, x_l)$.
The purpose of this note is to prove the following theorem,
which gives an answer of Poonen's question in
\cite{Po1}.

\begin{Theorem}
\label{thm:a:gen:lang:conj:fg:intro}
Let $\Gamma$ be a subgroup of finite rank in $A(\overline{K})$, and
$X$ a subvariety of $A_{\overline{K}}$.
Fix a basis $\{\gamma_1, \ldots, \gamma_n \}$ of $\Gamma \otimes \QQ$.
If the set $\{ x \in X(\overline{K}) \mid
\delta_L^{\overline{B}}(\gamma_1, \ldots, \gamma_n, x) \leq \epsilon \}$
is Zariski dense in $X$ for every positive number $\epsilon$,
then $X$ is a translation of an abelian subvariety of $A_{\overline{K}}$
by an element of $\Gamma_{div}$,
where $\Gamma_{div} = \{ x \in A(\overline{K}) \mid
\text{$nx \in \Gamma$ for some positive integer $n$} \}$.
\end{Theorem}

In the case where $d = 0$, Poonen proved the equivalent result in
\cite{Po1}.
Our argument for the proof of the above theorem essentially
follows his ideas. 
A new point is that we remove measure-theoretical arguments from his
original one, so that we can apply it to our case.
Finally, we note that
Theorem~\ref{thm:a:gen:lang:conj:fg:intro} substantially includes
Lang's conjecture in the absolute form:

\begin{LangConj}
Let $A$ be a complex abelian variety, 
$\Gamma$ a subgroup of finite rank in $A(\CC)$,
and $X$ a subvariety of $A$.
Then, there are abelian subvarieties $C_1, \ldots, C_n$
of $A$, and $\gamma_1, \ldots, \gamma_n \in \Gamma$
such that 
\[
\overline{X(\CC) \cap \Gamma} = \bigcup_{i=1}^n (C_i + \gamma_i)
\quad\text{and}\quad
X(\CC) \cap \Gamma = \bigcup_{i=1}^n (C_i(\CC) + \gamma_i) \cap \Gamma.
\]
\end{LangConj}

\renewcommand{\theTheorem}{\arabic{section}.\arabic{Theorem}}
\renewcommand{\theClaim}{\arabic{section}.\arabic{Theorem}.\arabic{Claim}}
\renewcommand{\theequation}{\arabic{section}.\arabic{Theorem}.\arabic{Claim}}
\section{Review of arithmetic height functions over finitely generated fields}
\label{sec:rev:arith:height:fun}
In this section, we give a quick review of arithmetic height functions
over finitely generated fields. For details, see \cite{MoArht}.

Let $K$ be a finitely generated field over $\QQ$ with
$d = \trdeg_{\QQ}(K)$, and let
$\overline{B} = (B; \overline{H}_1, \ldots, \overline{H}_d)$ 
be a big polarization of $K$, i.e.,
$B$ is a normal projective scheme over $\ZZ$, whose function
field is $K$, and $\overline{H}_1, \ldots, \overline{H}_d$
are nef and big $C^{\infty}$-hermitian line bundles on $B$.
For the definition of nef and big $C^{\infty}$-hermitian line bundles,
see \cite[\S2]{MoArht}.
Let $X$ be a projective variety over $K$ and $L$ a line bundle on $X$.
Let us consider a $C^{\infty}$-model $(\mathcal{X}, \overline{\mathcal{L}})$
of $(X, L)$ over $B$. Namely,
$\mathcal{X}$ is a projective integral scheme over $B$, whose generic fiber 
over $B$ is $X$,
and $\overline{\mathcal{L}}$ is a $C^{\infty}$-hermitian $\QQ$-line bundle on $\mathcal{X}$,
which gives rise to $L$ on the generic fiber of $\mathcal{X} \to B$.
For $x \in X(\overline{K})$, let $\Delta_x$ be the closure of the image
$\Spec(\overline{K}) \overset{x}{\longrightarrow} X \hookrightarrow \mathcal{X}$. Then,
we define the height of $x$ with respect to the polarization $\overline{B}$ and
the $C^{\infty}$-model $(\mathcal{X}, \overline{\mathcal{L}})$ to be
\[
h^{\overline{B}}_{(\mathcal{X}, \overline{\mathcal{L}})}(x) =
\frac{1}{[K(x):K]}\adeg \left(
\acherncl_1(\rest{\overline{\mathcal{L}}}{\Delta_x}) \cdot
\acherncl_1(\rest{\pi^* (\overline{H}_1)}{\Delta_x}) \cdots
\acherncl_1(\rest{\pi^* (\overline{H}_d)}{\Delta_x})
\right),
\]
where $\pi : \mathcal{X} \to B$ is the canonical morphism.
If $(\mathcal{X}', \overline{\mathcal{L}}')$ is another $C^{\infty}$-model of $(X, L)$,
then there is a constant $C$ such that
\[
\vert h^{\overline{B}}_{(\mathcal{X}, \overline{\mathcal{L}})}(x) -
h^{\overline{B}}_{(\mathcal{X}', \overline{\mathcal{L}}')}(x) \vert
\leq C
\]
for all $x \in X(\overline{K})$. Thus, modulo the set of bounded functions,
we can assign the unique height function $h^{\overline{B}}_L : X(\overline{K}) \to \RR$
to $\overline{B}$ and $L$.
Note that if $\sigma \in \Gal(\overline{K}/K)$, then
$\Delta_x = \Delta_{\sigma(x)}$. Thus, $h^{\overline{B}}_L(\sigma(x)) = h^{\overline{B}}_L(x)$.
The first important theorem is the following Northcott's theorem for our height
functions.

\begin{Theorem}[{\cite[Theorem~4.3]{MoArht}}]
\label{thm:northcott:thm:fun:field}
If $L$ is ample, then,
for any numbers $M$ and any positive integers $e$, the set
\[
\left\{ x \in X(\overline{K}) \mid h^{\overline{B}}_L(x) \leq M,
\quad [K(x) : K] \leq e \right\} 
\]
is finite.
\end{Theorem}

Let $A$ be an abelian variety over $K$, and $L$ a symmetric ample line bundle on $A$.
Then, as the usual height functions over a number field,
there is the canonical height function $\hat{h}^{\overline{B}}_L$.
This gives rise to a quadric form on $A(\overline{K})$, so that
if we set
\[
\langle x, y \rangle^{\overline{B}}_L =
\frac{1}{2} \left(
\hat{h}^{\overline{B}}_L(x + y) - \hat{h}^{\overline{B}}_L(x) - 
\hat{h}^{\overline{B}}_L(y)
\right)
\]
for $x, y \in A(\overline{K})$, then
$\langle \ , \ \rangle^{\overline{B}}_L$ is a bi-linear form on
$A(\overline{K})$.
Concerning this bi-linear form, we have the following.

\begin{Proposition}[{\cite[\S\S3.4]{MoArht}}]
\label{prop:element:height:abelian}
\begin{enumerate}
\renewcommand{\labelenumi}{(\arabic{enumi})}
\item
$\langle x, x \rangle^{\overline{B}}_L \geq 0$ for all $x \in A(\overline{K})$, and
the equality holds if and only if $x$ is a torsion point.
Namely, $\langle \ , \ \rangle^{\overline{B}}_L$ is positive definite
on $A(\overline{K}) \otimes \QQ$.

\item
If $f : A \to A'$ is a homomorphism of abelian varieties over $K$, and
$L'$ is a symmetric ample line bundle on $A'$, then there is a positive number
$a$ with
\[
 \langle f(x), f(x) \rangle^{\overline{B}}_{L'} \leq a \langle x, x \rangle^{\overline{B}}_L
\]
for all $x \in A(\overline{K})$.
\end{enumerate}
\end{Proposition}

\begin{Remark}
\label{rem:2:of:ineq:height}
(2) of Proposition~\ref{prop:element:height:abelian} holds even if
$f$, $A'$ and $L'$ are not defined over $K$. Let $K'$ be a finite extension field of
$K$ such that $f$, $A'$ and $L'$ are defined over $K'$.
Let $\phi : B^{K'} \to B$ be the normalization of $B$ in $K'$.
Then, $\overline{B}^{K'} = (B^{K'}; \phi^*(\overline{H}_1), \ldots, \phi^*(\overline{H}_d))$
gives rise to a big polarization of $K'$. Thus, there is a positive number $a'$
with
\[
 \langle f(x), f(x) \rangle^{\overline{B}^{K'}}_{L'} \leq 
a' \langle x, x \rangle^{\overline{B}^{K'}}_L
\]
for all $x \in A(\overline{K})$.
On the other hand, $\langle \ , \ \rangle^{\overline{B}^{K'}}_L = [K' : K]
\langle \ , \ \rangle^{\overline{B}}_L$. Hence,
\[
 \langle f(x), f(x) \rangle^{\overline{B}^{K'}}_{L'} \leq 
a' [K' : K] \langle x, x \rangle^{\overline{B}}_L
\]
for all $x \in A(\overline{K})$.
\end{Remark}

The crucial result for this note is the following solution of Bogomolov's conjecture
over finitely generated fields, which is a generalization of
\cite{UlPos} and \cite{ZhEqui}.

\begin{Theorem}[{\cite[Theorem~8.1]{MoArht}}]
\label{thm:bogomolov:conj:fun}
Let $X$ be a subvariety of $A_{\overline{K}}$.
If the set
\[
 \{ x \in X(\overline{K}) \mid \hat{h}^{\overline{B}}_L(x) \leq \epsilon \}
\]
is Zariski dense in $X$ for every positive number $\epsilon$,
then $X$ is a translation of an abelian subvariety of $A_{\overline{K}}$ 
by a torsion point.
\end{Theorem}

\section{Small points with respect to a group of finite rank}
\label{sec:small:points}
The contexts in this section are essentially due to Poonen \cite{Po1}.
We just deal with his ideas in a general situation.

Let $K$ be a finitely generated field over $\QQ$ with
$d = \trdeg_{\QQ}(K)$, and let
$\overline{B}$ 
be a big polarization of $K$.
Let $A$ be an abelian variety over $K$, and $L$ a symmetric
ample line bundle on $A$. Let
\[
\langle\ , \ \rangle_L^{\overline{B}} :
A(\overline{K}) \times A(\overline{K}) \to \RR
\]
be the height pairing associated with $\overline{B}$ and $L$
as in \S\ref{sec:rev:arith:height:fun}.
Let $\Gamma$ be a subgroup of finite rank in $A(\overline{K})$.
A non-empty subset $S$ of $A(\overline{K})$ is said to
be {\em small with respect to $\Gamma$}
if there is a decomposition
$s = \gamma(s) + z(s)$ for each $s \in S$ with the following properties:
\begin{enumerate}
\renewcommand{\labelenumi}{(\alph{enumi})}
\item
$\gamma(s) \in \Gamma$ for all $s \in S$.

\item
For any $\epsilon > 0$,
there is a finite proper subset $S'$ of $S$ such that
$\langle z(s), z(s) \rangle_{L}^{\overline{B}} \leq \epsilon$
for all $s \in S \setminus S'$.
\end{enumerate}
Especially a small subset $S$ with respect to $\{ 0 \}$
is said to be {\em small}.
Namely, 
a non-empty subset $S$ of $A(\overline{K})$ is small if and only if,
for any positive numbers $\epsilon$, there is a finite proper subset
$S'$ of $S$ with $\langle x, x \rangle_L^{\overline{B}} \leq \epsilon$
for all $s \in S \setminus S'$. 
Note that in the above definition, $S'$ is proper, i.e.,
$S \setminus S' \not= \emptyset$.
Let us begin with the following proposition.

\begin{Proposition}
\label{prop:element:prop:small}
Let $S$ be a non-empty subset of $A(\overline{K})$ and $\Gamma$ a subgroup of finite rank
in $A(\overline{K})$.
Then, we have the following:
\begin{enumerate}
\renewcommand{\labelenumi}{(\arabic{enumi})}
\item
If $S$ is small with respect to $\Gamma$, then
any infinite subsets of $S$ are small with respect to $\Gamma$.

\item
We assume that $S$ is finite. Then $S$ is small
\rom{(}with respect to $\{ 0 \}$\rom{)} if and only if
$S$ contains a torsion point.

\item
We assume that $S$ is infinite.
Let $N$ be a positive integer, and
$[N]$ an endomorphism of $A$ given by $[N](x) = Nx$.
If $S$ is small with respect to $\Gamma$, then
so is $[N](S)$.

\item
Let $\{ x_n \}$ be a sequence in $A(\overline{K})$ with the following properties:
\begin{enumerate}
\renewcommand{\labelenumii}{(\arabic{enumi}.\arabic{enumii})}
\item
If $n \not= m$, then $x_n \not= x_m$.

\item
Each $x_n$ has a decomposition $x_n = \gamma_n + y_n$ with $\gamma_n \in \Gamma$.

\item
${\displaystyle \lim_{n \to \infty} \langle y_n, y_n \rangle_L^{\overline{B}} = 0}$.
\end{enumerate}
Then, $\{ x_n \mid n=1, 2, \ldots \}$ is small with respect to $\Gamma$.
\end{enumerate}
\end{Proposition}

\Proof
(1) and (4) are obvious.

\medskip
(2) Clearly, if $S$ contains a torsion point, then $S$ is small.
We assume that $S$ is small. We set
$\lambda = \min \{ \langle s, s \rangle_{L}^{\overline{B}} \mid s \in S \}$.
If $\lambda > 0$, then there is a finite proper subset $S'$ of $S$
such that
$\langle s, s \rangle_{L}^{\overline{B}} < \lambda$ for all $s \in S \setminus S'$.
This is a contradiction. Thus, $\lambda = 0$, which means that
$S$ contains a torsion point.

\medskip
(3) We fix a map $t : [N](S) \to S$ with $[N](t(s)) = s$ for all $s \in [N](S)$.
Then, we have a decomposition $s = [N](\gamma(t(s))) + [N](z(t(s)))$ for each $s \in [N](S)$.
Clearly (a) in the definition of small sets is satisfied.
Let $\epsilon$ be an arbitrary positive number. Then,
there is a finite subset $T$ of $S$ such that
$\langle z(s), z(s) \rangle_{L}^{\overline{B}} \leq \epsilon/N^2$
for all $s \in S \setminus T$.
If we set $T' = \{ s \in [N](S) \mid t(s) \in T \}$, then $T'$ is finite.
Moreover, 
$\left\langle [N](z(t(s))), [N](z(t(s))) \right\rangle_{L}^{\overline{B}} \leq \epsilon$
for all $s \in [N](S) \setminus T'$.
Therefore, we have (b) in the definition of small sets.
\QED

Moreover, we have the following, which is a consequence of Bogomolov's conjecture.

\begin{Theorem}
\label{thm:conq:bogomolov:conj}
Let $S$ be a small set of $A(\overline{K})$, i.e.,
$S$ is small with respect to $\{ 0 \}$.
Then, there are abelian subvarieties $C_1, \ldots, C_n$,
torsion points $c_1, \ldots, c_n$, and finite non-torsion points $b_1, \ldots, b_m$
such that
\[
\overline{S} = \bigcup_{i=1}^n (C_i + c_i) \cup \{ b_1, \ldots, b_m \},
\]
where $\overline{S}$ is the Zariski closure of $S$.
\end{Theorem}

\Proof
It is sufficient to show that
a positive dimensional irreducible component $X$ of $\overline{S}$ is
a translation of an abelian subvariety of $A$ by a torsion point.
Let $S'$ be the set of points in $S$, which is contained in $X(\overline{K})$.
Then, the Zariski closure of $S'$ is $X$. In particular, $S'$ is infinite set,
so that $S'$ is small. Thus, $X$ is a translation of an abelian subvariety of $A$ by a torsion point
by virtue of Theorem~\ref{thm:bogomolov:conj:fun}.
\QED

Let $S$ be a small subset with respect to $\Gamma$.
For each $n \geq 2$, let us consider a homomorphism $\beta_n : A^n \to A^{n-1}$
given by $\beta_n(a_1, \ldots a_n) = (a_2 - a_1, a_3 - a_1, \ldots, a_n - a_1)$.
Let $F$ be a finite extension field of $K$ in $\overline{K}$.
For $x \in A(\overline{K})$, we denote by $O_{F}(x)$ the orbit of $x$ by the Galois group
$\Gal(\overline{K}/F)$. Noting $O_{F}(x)^n \subseteq A(\overline{K})^n$,
for a subset $T$ of $S$,
we define $\DD_n(T, F)$ to be
\[
 \DD_n(T, F) = \bigcup_{s \in T} \beta_n(O_{F}(s)^n).
\]
We denote the Zariski closure of $\DD_n(T, F)$ by $\ZDD_n(T, F)$.
On $A^n$, we can give the height pairing associated with
$\bigotimes_{i=1}^n p_i^*(L)$ and $\overline{B}$, where $p_i : A^n \to A$
is the projection to the $i$-th factor.
By abuse of notation, we denote this by $\langle\ , \ \rangle_L^{\overline{B}}$. 

\begin{Proposition}
\label{prop:small:d:n}
Let $f : A \to A'$ be a homomorphism of abelian varieties over $\overline{K}$.
Let $F$ be a finite extension field of $K$ in $\overline{K}$.
We assume that there is a finitely generated subgroup $\Gamma_0$ of $\Gamma$
such that $\Gamma_0 \subseteq A(K)$ and $\Gamma_0 \otimes \QQ = \Gamma \otimes \QQ$.
Then, we have the following:
\begin{enumerate}
\renewcommand{\labelenumi}{(\arabic{enumi})}
\item
$f^{n-1}(\DD_n(S, F))$ is small \rom{(}with respect to $\{ 0 \}$\rom{)}, where
$f^{n-1} : A^{n-1} \to {A'}^{n-1}$ is the morphism given by 
$f^{n-1}(x_1, \ldots, x_{n-1}) = (f(x_1), \ldots, f(x_{n-1}))$.

\item
Let $b_1, \ldots, b_l$ be non-torsion points in $f^{n-1}(\DD_n(S, F))$.
Then, there is a finite proper subset $S'$ of $S$ such that
$b_i \not\in f^{n-1}(\DD_n(S \setminus S', F))$ for all $i$.
\end{enumerate}
\end{Proposition}

\Proof
Let $\sigma, \tau$ be elements of $\Gal(\overline{K}/F)$.
Then, $\sigma(\gamma(s)) - \tau(\gamma(s))$ is torsion
because $n \gamma(s) \in \Gamma_0$
for some $n > 0$.
Thus,
\[
\Vert \sigma(s) - \tau(s) \Vert_L^{\overline{B}} =
\Vert \sigma(z(s)) - \tau(z(s)) \Vert_L^{\overline{B}} 
\leq 2\Vert z(s) \Vert_L^{\overline{B}},
\]
where $\Vert x \Vert_L^{\overline{B}} = \sqrt{\langle x,x \rangle_L^{\overline{B}}}$.
Therefore,
\[
\Vert \beta_n(x) \Vert_L^{\overline{B}} \leq 2 \sqrt{n-1} \Vert z(s) \Vert_L^{\overline{B}}
\]
for all $x \in O_{F}(s)^n$.
Let $L'$ be a symmetric ample line bundle on $A'$.
Then, by (2) of Proposition~\ref{prop:element:height:abelian}
(or Remark~\ref{rem:2:of:ineq:height}),
there is a positive constant $a$ with
$\langle f(x), f(x) \rangle_{L'}^{\overline{B}} \leq a 
\langle x, x \rangle_{L}^{\overline{B}}$ for all $x \in A(\overline{K})$.
Thus,
\addtocounter{Claim}{1}
\begin{equation}
\label{prop:small:d:n:eqn:1}
\Vert f^{n-1}(\beta_n(x)) \Vert_{L'}^{\overline{B}} 
\leq 2\sqrt{a (n-1)} \Vert z(s) \Vert_L^{\overline{B}}
\end{equation}
for all $x \in O_{F}(s)^n$.

First, let us see (2).
We set $\mu = \min \{ \Vert b_i \Vert_{L'}^{\overline{B}} \mid i =1, \ldots, l \} > 0$.
Then there is a finite proper subset $S'$ of $S$ with
\[
\Vert z(s) \Vert_L^{\overline{B}} < \frac{\mu}{2\sqrt{a(n-1)}}
\]
for all $s \in S \setminus S'$.
Thus, by \eqref{prop:small:d:n:eqn:1},
\[
\Vert f^{n-1}(\beta_n(x)) \Vert_{L'}^{\overline{B}} < \mu
\]
for all $x \in \bigcup_{s \in S \setminus S'} O_{F}(s)^n$.
Hence, $b_i \not\in f^{n-1}(\DD_n(S \setminus S', F))$
for all $i$.

Next we consider (1).
If $f^{n-1}(\DD_n(S, F))$ is infinite, then the assertion of (1) is obvious
by \eqref{prop:small:d:n:eqn:1}.
Otherwise, let $\{ b_1, \ldots, b_n \}$ be the set of
all non-torsion points in $f^{n-1}(\DD_n(S, F))$.
Then, by (2), we can find a finite proper subset $S'$ of $S$ with
\[
\emptyset \not= f^{n-1}(\DD_n(S \setminus S', F)) \subseteq 
f^{n-1}(\DD_n(S, F)) \setminus \{ b_1, \ldots, b_n \}.
\]
Hence $f^{n-1}(\DD_n(S, F))$ contains a torsion point.
Therefore, $f^{n-1}(\DD_n(S, F))$ is small.
\QED

Let $S$ be a small subset with respect to $\Gamma$.
From now on, we assume the following:
\begin{enumerate}
\renewcommand{\labelenumi}{(\Alph{enumi})}
\item
$S$ is infinite.

\item
There is a finitely generated subgroup $\Gamma_0$ of $\Gamma$
such that $\Gamma_0 \subseteq A(K)$ and $\Gamma_0 \otimes \QQ = \Gamma \otimes \QQ$.
\end{enumerate}
Let $F$ be a finite extension field of $K$ in $\overline{K}$.
A pair $(S, F)$ is said to be {\em $n$-minimized} if
the following properties are satisfied:
\begin{enumerate}
\renewcommand{\labelenumi}{(\roman{enumi})}
\item
$\ZDD_n(S', F') = \ZDD_n(S, F)$
for any infinite subsets $S'$ of $S$ and any finite extension fields $F'$ of $F$
in $\overline{K}$. (Recall that
$\ZDD_n(\cdot, \cdot)$ is the Zariski closure of $\DD_n(\cdot, \cdot)$.)

\item
$\ZDD_n([N](S), F) = \ZDD_n(S, F)$
for any positive integers $N$.
\end{enumerate}
Note that $[N](O_F(s)) = O_F([N](s))$ for $s \in S$ and a positive integer $N$,
so that $\DD_n([N](S), F) = [N](\DD_n(S, F))$.
Therefore, (ii) is equivalent to saying that
$[N](\ZDD_n(S, F)) = \ZDD_n(S, F)$
for any positive integers $N$.
First let us consider the following proposition.

\begin{Proposition}
\label{prop:property:mini}
\begin{enumerate}
\renewcommand{\labelenumi}{(\arabic{enumi})}
\item
If we fix $n \geq 2$, then there are an infinite subset $T$ of $S$, a positive integer $N$, and
a finite extension field $F$ 
of $K$ in $\overline{K}$ such that $([N](T), F)$ is $n$-minimized.

\item
Let $F$ be a finite extension field of $K$ in $\overline{K}$.
Let $N$ be a positive integer, $S'$ an infinite subset of $[N](S)$,
and $F'$ a finite extension field of $F$ in $\overline{K}$.
If $(S, F)$ is $n$-minimized, then
$\ZDD_n(S', F') = \ZDD_n(S, F)$.
\end{enumerate}
\end{Proposition}

\Proof
(1)
Let $F$ be a finite extension field of $K$ in $\overline{K}$.
A pair $(S, F)$ is said to be {\em weakly $n$-minimized} if
the above property (i) is satisfied.
First, we claim the following.

\begin{Claim}
\label{claim:prop:property:mini}
\begin{enumerate}
\renewcommand{\labelenumi}{(\alph{enumi})}
\item
If we fix $n \geq 2$, then there are an infinite subset $T$ of $S$ and a finite extension field $F$
of $K$ such that $(T, F)$ is weakly $n$-minimized.

\item
Let $F$ be a finite extension field of $K$ in $\overline{K}$.
If $(S, F)$ is weakly $n$-minimized, then
there are abelian subvarieties $C_1, \ldots, C_n$, and
torsion points $c_1, \ldots, c_n$ such that
\[
\ZDD_n(S, F) = \bigcup_{i=1}^n (C_i + c_i).
\]

\item
Let $F$ be a finite extension field of $K$ in $\overline{K}$,
and $N$ a positive integer.
If $(S, F)$ is weakly $n$-minimized, then
so is $([N](S), F)$.
\end{enumerate}
\end{Claim}

(a) This is obvious by Noetherian induction.

\medskip
(b) By Theorem~\ref{thm:conq:bogomolov:conj},
there are abelian subvarieties $C_1, \ldots, C_n$,
torsion points $c_1, \ldots, c_n$, and finite non-torsion
points $b_1, \ldots, b_m$
such that
\[
\ZDD_n(S, F) = \bigcup_{i=1}^n (C_i + c_i) \cup \{ b_1, \ldots, b_m \}.
\]
By virtue of (2) of Proposition~\ref{prop:small:d:n},
we can find a finite set $T$ of $S$ such that
\[
\ZDD_n(S \setminus T, K) \subseteq \bigcup_{i=1}^n (C_i + c_i)
\subseteq \ZDD_n(S, F).
\]
Here, $\ZDD_n(S \setminus T, K) = \ZDD_n(S, K)$.
Thus, we get (b).

\medskip
(c) Let $S_1$ be an infinite subset of $[N](S)$ and $F'$ a finite extension field of $F$
in $\overline{K}$.
We take a subset $S'$ of $S$ with $[N](S') = S_1$.
Then, $\ZDD_n(S', F') = \ZDD_n(S, F)$.
Thus, since $[N]$ is a finite and surjective morphism, we can see
\[
\ZDD_n(S_1, F') = \ZDD_n([N](S'), F') =
[N](\ZDD_n(S', F')) = [N](\ZDD_n(S, F))
= \ZDD_n([N](S), F).
\]
Hence, we have (c).

\bigskip
Let us start the proof of (1).
By virtue of (a), there are an infinite subset $T$ of $S$ and a finite extension field $F$
of $K$ such that $(T, F)$ is weakly $n$-minimized.
Hence, by (b),
there are abelian subvarieties $C_1, \ldots, C_n$, and
torsion points $c_1, \ldots, c_n$ such that
\[
\ZDD_n(T, F) = \bigcup_{i=1}^n (C_i + c_i).
\]
Let $N$ be a positive integer with $Nc_i = 0$ for all $i$.
Then,
\[
\ZDD_n([N](T), F) = [N](\ZDD_n(T, F)) =
\bigcup_{i=1}^n C_i.
\]
Here we claim that $([N](T), F)$ is $n$-minimized.
By (c), $([N](T), F)$ is weakly $n$-minimized.
Moreover, for any positive integers $N'$,
\begin{align*}
\ZDD_n([N']([N](T)), F) & = [N'](\ZDD_n([N](T), F)) \\
& =[N']\left(\bigcup_{i=1}^n C_i\right) = \bigcup_{i=1}^n C_i \\
& = \ZDD_n([N](T), F).
\end{align*}
Thus, $([N](T), F)$ is $n$-minimized.

\bigskip
(2)
Let $N$ be a positive integer, $S'$ an infinite subset of $[N](S)$, and
$F'$ a finite extension field of $F$.
By (c), $([N](S), F)$ is weakly $n$-minimized.
Thus,
\[
\ZDD_n(S', F') = \ZDD_n([N](S), F) = \ZDD_n(S, F).
\]
Therefore, we get (2).
\QED

Finally, let us consider the following theorem, which is crucial
for our note.

\begin{Theorem}
\label{thm:equiv:cond:n:min}
Let $F$ be a finite extension field of $K$ in $\overline{K}$.
Then, the following \rom{(1)}, \rom{(2)} and \rom{(3)} are equivalent.
\begin{enumerate}
\renewcommand{\labelenumi}{(\arabic{enumi})}
\item
$(S, F)$ is $n$-minimized for all $n \geq 2$.

\item
$(S, F)$ is $n$-minimized for some $n \geq 2$.

\item
$(S, F)$ is $2$-minimized.
\end{enumerate}
Moreover, under the above equivalent conditions,
there is an abelian subvariety $C$ of $A_{\overline{K}}$
such that $\ZDD_n(S, F) = C^{n-1}$ for all $n \geq 2$.
\end{Theorem}

\Proof
Let us begin with the following two lemmas.

\begin{Lemma}
\label{lem:for:prop:2:n:minimal}
Let $F$ be a finite extension field of $K$ in $\overline{K}$, and
$C$ an abelian subscheme of $A_F$ over $F$.
We assume that there is a positive integer $e$ with the following property:
For each $s \in S$, there is a subset $T(s)$ of $O_F(s) \times O_F(s)$
such that $\beta_2(T(s)) \subseteq C(\overline{K})$ and
$\#(T(s)) \geq \#(O_F(s) \times O_F(s))/e$.
Then, there is a finite subset $S'$ of $S$ and
a positive integer $N$ with
$\DD_2([N](S \setminus S'), F) \subseteq C(\overline{K})$.
\end{Lemma}

\Proof
Let $\pi : A \to A/C$ be a natural homomorphism.
Fix $s \in S$.
Let $F'$ be a finite Galois extension of $F$ such that
$F'$ contains $F(s)$.
Then, there is a natural surjective map
\[
\phi : \Gal(F'/F) \to O_F(s),
\]
whose fibers are cosets of the stabilizer of $s$.
If we set $E(s) = (\phi \times \phi)^{-1}(T(s))$, then
$\#(E(s)) \geq \#(\Gal(F'/F) \times \Gal(F'/F))/e$ and
$\sigma(\pi(s)) = \tau(\pi(s))$ for all $(\sigma, \tau) \in E(s)$.
Let $G_{\pi(s)}$ be the stabilizer of $\pi(s)$ by the action of $\Gal(F'/F)$, and
let $R$ be the set of all $(\sigma, \tau) \in \Gal(F'/F) \times \Gal(F'/F)$
with $\sigma(\pi(s)) = \tau(\pi(s))$.
Then, we have
\[
\#(R) = \#(G_{\pi(s)}) \#(\Gal(F'/F))
\quad\text{and}\quad
\#(R) \geq  \frac{\#(\Gal(F'/F) \times \Gal(F'/F))}{e}.
\]
Thus, $[\Gal(F'/F):G_{\pi(s)}] \leq e$, which means that
$[F(\pi(s)):F] \leq e$.
Then, since $\pi(\DD_2(S, F))$ is small,
by virtue of Northcott's theorem (cf. Theorem~\ref{thm:northcott:thm:fun:field}), 
$\pi(\DD_2(S, F))$ is finite.
By (2) of Proposition~\ref{prop:small:d:n},
there is a finite proper subset $S'$ of $S$ such that
$\pi(\DD_2(S \setminus S', F))$ consists of torsion points.
Hence, there is a positive integer $N$ such that
$[N](\pi(\DD_2(S \setminus S', F))) = \{ 0 \}$.
Therefore, $\DD_2([N](S \setminus S'), F) \subseteq C(\overline{K})$.
\QED

\begin{Lemma}
\label{lem:structure:2:n:minimal}
Let $F$ be a finite extension field of $K$ in $\overline{K}$.
If $(S, F)$ are $2$-minimized,
then there is an abelian subvariety $C$ of $A_{\overline{K}}$
such that $\ZDD_n(S, F) = C^{n-1}$ for all $n \geq 2$.
\end{Lemma}

\Proof
First, let us consider the case $n=2$.
By using (b) of Claim~\ref{claim:prop:property:mini},
we can find abelian subvarieties $C_1, \ldots, C_e$ with
\[
\ZDD_2(S, F) = \bigcup_{i=1}^e C_i
\]
because $\ZDD_2(S, F)$ is stable by the endomorphism
$[N]$ for every positive integer $N$.
Thus, in order to see $e = 1$, it is sufficient to find $C_i$,
a positive integer $N_1$, an infinite subset $S_1$ of $S$, 
and a finite extension field $F_1$ of $F$
such that
\[
\DD_2([N_1](S_1), F_1) \subseteq C_i(\overline{K}).
\]
Let $F_1$ be a finite extension field of $F$ such that $C_i$'s are defined over $F_1$.
For each $s \in S$, let $T_i(s)$ be the set of all elements $x \in O_{F_1}(s)^2$
with $\beta_2(x) \in C_i(\overline{K})$.
We choose a map $\lambda : S \to \{1, \ldots, e \}$ such that
$\#(T_{\lambda(s)}(s))$ gives rise to the maximal value
in $\{ \#(T_{i}(s)) \mid i=1, \ldots, e \}$.
By using the pigeonhole principle, 
there are $i \in \{1, \ldots, e \}$ and an infinite subset $S'$ of $S$
with $\lambda(s) = i$ for all $s \in S'$.
Then, for all $s \in S'$,
$\beta_2(T_i(s)) \subseteq C_i(\overline{K})$ and
$\#(T_i(s)) \geq \#(O_{F_1}(s)^2)/e$.
Thus, by Lemma~\ref{lem:for:prop:2:n:minimal},
there are an infinite subset $S_1$ of $S'$ and a positive integer $N_1$
with $\DD_2([N_1](S_1), F_1) \subseteq C_i(\overline{K})$.

\medskip
From now on, we denote $C_i$ by $C$. Then, $\ZDD_2(S, F) = C$.
Let us try to see $\ZDD_n(S, F) = C^{n-1}$ for all $n \geq 2$.
Clearly, $\ZDD_n(S, F) \subseteq C^{n-1}$.
Thus it is sufficient to find
a positive integer $N_2$, an infinite subset $S_2$ of $S$, and a finite extension field $F_2$ of $F$
such that
\[
\ZDD_n([N_2](S_2), F_2) = C^{n-1}.
\]
By (1) of Proposition~\ref{prop:property:mini},
there are a positive integer $N_2$, an infinite subset $S_2$ of $S$ and 
a finite extension field $F_2$ of $F$ such that
$([N_2](S_2), F_2)$ is $n$-minimized. 
Thus, as before,
there are abelian subvarieties $G_1, \ldots, G_l$ with
$\ZDD_n([N_2](S_2), F_2) = \bigcup_{j=1}^l G_j$.
Moreover, replacing $F_2$ by a finite extension field of $F_2$,
we may assume that $C$ and $G_j$'s are defined over $F_2$.
On this stage, we would like to show that
\[
\ZDD_n([N_2](S_2), F_2) = C^{n-1}.
\]

In the same way as before,
we can find $G_j$, say $G$, and
an infinite subset $S'$ of $[N_2](S_2)$ such that
for all $s \in S'$, there is a subset $T(s)$ of $O_{F_2}(s)^n$ with
$\#(T(s)) \geq \#(O_{F_2}(s)^n)/l$ and $\beta_n(T(s)) \subseteq G(\overline{K})$.
Let $C^{(q)} = 0 \times \cdots \times C \times \cdots \times 0$ be
the $q$-th factor of $C^{n-1}$, and
$G^{(q)} = G \cap C^{(q)}$ for $1 \leq q \leq n-1$. 
Since $G \subseteq C^{n-1}$, it is sufficient to see 
the following claim to conclude the proof of our lemma.

\begin{Claim}
$G^{(q)} = C^{(q)}$ for each $1 \leq q \leq n-1$.
\end{Claim}

For each $t_1, \ldots, t_{q}, t_{q+2}, \ldots, t_{n} \in O_{F_2}(s)$,
we set
\[
J(t_1, \ldots, t_{q}, t_{q+2}, \ldots, t_{n}) =
\{ x \in O_{F_2}(s) \mid  (t_1, \ldots, t_{q}, x, t_{q+2}, \ldots, t_{n}) \in T(s) \}.
\]
We choose $s_1, \ldots, s_{q}, s_{q+2}, \ldots, s_{n} \in O_{F_2}(s)$ such that
$\#(J(s_1, \ldots, s_{q}, s_{q+2}, \ldots, s_{n}))$ is maximal among
$\{ \#(J(t_1, \ldots, t_{q}, t_{q+2}, \ldots, t_{n})) \mid 
t_1, \ldots, t_{q}, t_{q+2}, \ldots, t_{n} \in O_{F_2}(s) \}$. Then,
\[
\#(J(s_1, \ldots, s_{q}, s_{q+2}, \ldots, s_{n})) \#( O_{F_2}(s)^{n-1} )
\geq \#(T(s)) \geq \frac{\#(O_{F_2}(s)^n)}{l}.
\]
Thus if we set $L(s) = J(s_1, \ldots, s_{q}, s_{q+2}, \ldots, s_{n})$, then
$\#(L(s)) \geq \#(O_{F_2}(s))/l$ and
\[
\beta_n(s_1, \ldots, s_{q}, x, s_{q+2}, \ldots, s_{n}) \in G(\overline{K})
\]
for all $x \in L(s)$.
Therefore, for all $(x, x') \in L(s) \times L(s)$,
\begin{multline*}
\beta_n(0, \ldots, 0, x - x', 0, \ldots, 0) = \\
\beta_n(s_1, \ldots, s_{q}, x, s_{q+2}, \ldots, s_{n}) -
\beta_n(s_1, \ldots, s_{q}, x, s_{q+2}, \ldots, s_{n}) \in  G(\overline{K}).
\end{multline*}
This means that
$\beta_2(x,  x') \in G^{(q)}(\overline{K})$ for all
$(x, x') \in L(s) \times L(s)$ if we view $G^{(q)}$ as a subscheme of $A$.
Here $\#(L(s) \times L(s)) \geq \#(O_{F_2}(s) \times O_{F_2}(s))/l^2$.
By Lemma~\ref{lem:for:prop:2:n:minimal},
there are an infinite subset $S''$ of $S'$ and a positive integer $N''$
with $\ZDD_2([N''](S''), F_2) \subseteq G^{(q)}$, which implies that $G^{(q)} = C^{(q)}$
because $\ZDD_2([N''](S''), F_2) = C$
by (2) of Proposition~\ref{prop:property:mini}.
\QED

\bigskip
Let us start the proof of Theorem~\ref{thm:equiv:cond:n:min}.
The last assertion is nothing more than
Lemma~\ref{lem:structure:2:n:minimal}, 
so that it is sufficient to show that (2) $\Longrightarrow$ (3) and
(3) $\Longrightarrow$ (1).

\medskip
(2) $\Longrightarrow$ (3):
By (1) of Proposition~\ref{prop:property:mini},
there are an infinite subset $T$ of $S$, a positive integer $N_1$, and
a finite extension field $F_1$ of $F$ 
in $\overline{K}$ such that $([N_1](T), F_1)$ is $2$-minimized.
Then, by Lemma~\ref{lem:structure:2:n:minimal},
there is an abelian subvariety $C$ of $A_{\overline{K}}$
such that $\ZDD_2([N_1](T), F_1) = C$ 
and $\ZDD_n([N_1](T), F_1) = C^{n-1}$.
Thus, $\ZDD_n(S, F) = C^{n-1}$ because
$(S, F)$ is $n$-minimized.
For all $x, x' \in O_F(s)$ with $s \in S$,
\[
\beta_n(s, x, s, \ldots, s) - \beta_n(s, x', s, \ldots, s) =
(x - x', 0, \ldots, 0) \in C(\overline{K})^{n-1}.
\]
Thus, $\beta_2(O_{F}(s)^2) \subseteq C(\overline{K})$ for all $s \in S$. 
Therefore, $\ZDD_2(S, F) \subseteq C$.
Let $S'$ be an infinite subset of $S$, and
$F'$ a finite extension field of $K$. In order to see that
$\ZDD_2(S', F') = C$, we may assume that
$S' \subseteq T$ and $F_1 \subseteq F'$.
Then,
\[
[N_1](\ZDD_2(S', F')) = 
\ZDD_2([N_1](S'), F') = \ZDD_2([N_1](T), F_1) = C.
\]
Thus, $\ZDD_2(S', F') = C$
because $\ZDD_2(S', F') \subseteq C$.
Hence $(S, F)$ satisfies the property (i) in
the definition of ``2-minimized''.
Moreover, $[N](\ZDD_2(S, F)) = [N](C) = C$
for all positive integers $N$.
Therefore, $(S, F)$ is $2$-minimized.

\medskip
(3) $\Longrightarrow$ (1):
By Lemma~\ref{lem:structure:2:n:minimal},
there is an abelian subvariety $C$ of $A_{\overline{K}}$
such that $\ZDD_n(S, F) = C^{n-1}$ for all $n \geq 2$.
Fix $n \geq 2$.
By (1) of Proposition~\ref{prop:property:mini},
there are an infinite subset $T$ of $S$, a positive integer $N_1$, and
a finite extension field $F_1$ of $F$ 
in $\overline{K}$ such that $([N_1](T), F_1)$ is $n$-minimized.
Since $([N_1](T), F_1)$ is $2$-minimized and $\ZDD_2([N_1](T), F_1) = C$,
we have $\ZDD_n([N_1](T), F_1) = C^{n-1}$
by Lemma~\ref{lem:structure:2:n:minimal}.
Thus, as before, we can see that
$(S, F)$ is $n$-minimized.
\QED

\renewcommand{\theTheorem}{\arabic{section}.\arabic{subsection}.\arabic{Theorem}}
\renewcommand{\theClaim}{\arabic{section}.\arabic{subsection}.\arabic{Theorem}.\arabic{Claim}}
\renewcommand{\theequation}{\arabic{section}.\arabic{subsection}.\arabic{Theorem}.\arabic{Claim}}
\section{Proof of Theorem~\ref{thm:a:gen:lang:conj:fg:intro}}

\setcounter{Theorem}{0}
\subsection{Preliminary of linear algebra}
\label{subsec:pre:linear:algebra}
Let $V$ be a vector space over $\RR$, and
$\langle\ , \ \rangle$ an inner product on $V$.
For a finite set of linearly independent
vectors $\Lambda = \{ v_1, \ldots, v_n \}$,
we define
\[
\Delta_{\Lambda} : V \times V \to \RR
\]
to be
\[
\Delta_{\Lambda}(x, y) = \det
\begin{pmatrix}
\langle v_1, v_1 \rangle & \cdots & \langle v_1, v_n \rangle & \langle v_1, y \rangle \\
\vdots   & \ddots & \vdots & \vdots \\
\langle v_n, v_1 \rangle & \cdots & \langle v_n, v_n \rangle & \langle v_n, y \rangle \\
\langle x, v_1 \rangle & \cdots & \langle x, v_n \rangle & \langle x, y \rangle
\end{pmatrix}
\]
Then, we have the following:

\begin{Proposition}
\label{prop:elem:prop:B:lambda}
\begin{enumerate}
\renewcommand{\labelenumi}{(\arabic{enumi})}
\item
$\Delta_{\Lambda}$ is a bi-linear map.

\item
$\Delta_{\Lambda}$ is symmetric and positive semidefinite.

\item
For all $v \in \Span(\Lambda)$ and $x \in V$,
$\Delta_{\Lambda}(v, x) = 0$.

\item
If $\Lambda' = \{ v'_1, \ldots, v'_n \}$ is another finite
set of linearly independent vectors with
$\Span(\Lambda') = \Span(\Lambda)$, then
\[
\Delta_{\Lambda'} = \frac{\det(\Lambda')}{\det(\Lambda)} \Delta_{\Lambda},
\]
where $\det(\Lambda) = \det(\langle v_i, v_j \rangle)$ and
$\det(\Lambda') = \det(\langle v'_i, v'_j \rangle)$.

\item
There are linear maps $p_{\Lambda} : V \to \Span(\Lambda)$ and
$q_{\Lambda} : V \to \Span(\Lambda)^{\perp}$ with
$x =  p_{\Lambda}(x) + q_{\Lambda}(x)$ for all $x \in V$,
where $\Span(\Lambda)^{\perp} = \{ x \in V \mid \text{$\langle x, v \rangle = 0$
for all $v \in \Span(\Lambda)$} \}$.

\item
$\Delta_{\Lambda}(x, x) =
\det (\Lambda) \langle q_{\Lambda}(x), q_{\Lambda}(x) \rangle$
for all $x \in V$.
In particular,
$\Delta_{\Lambda}(x, x) \leq \det(\Lambda) \langle x, x \rangle$
and the equality holds if and only if
$x \in \Span(\Lambda)^{\perp}$.

\end{enumerate}
\end{Proposition}

\Proof
(1), (2) and (3) are straightforward from the definition of $\Delta_{\Lambda}$.

(4)
First of all, there is an invertible matrix $P$ with
$(v'_1, \ldots, v'_n) = (v_1, \ldots, v_n)P$. Then it is easy to see that
$(\langle v'_i, v'_j \rangle) = P (\langle v_i, v_j \rangle){}^tP$.
Thus, $\det(\Lambda') = \det(P)^2 \det(\Lambda)$.
On the other hand, since
\[
(v'_1, \ldots, v'_n, x) = (v_1, \ldots, v_n, x)
\begin{pmatrix} P & 0 \\ 0 & 1 \end{pmatrix},
\]
in the same way as above, we have
$\Delta_{\Lambda'}(x, x) = \det(P)^2 \Delta_{\Lambda}(x, x)$.
Therefore,
\[
\Delta_{\Lambda'}(x, y) =  \frac{\det(\Lambda')}{\det(\Lambda)} \Delta_{\Lambda}(x, y)
\]
for all $x, y \in V$ because
$2 \Delta_{\Lambda}(x, y) = \Delta_{\Lambda}(x+y, x+y) - \Delta_{\Lambda}(x, x) - 
\Delta_{\Lambda}(y, y)$.

(5)
For $x \in V$, solving the equation
\[
\sum_{j=1}^n \lambda_j \langle v_j, v_i \rangle = \langle x, v_i \rangle
\quad\text{for all $i =1, \ldots, n$},
\]
we can find a unique vector $v = \sum \lambda_j v_j \in \Span(\Lambda)$ such that
$x - v$ is perpendicular to $\Span(\Lambda)$.
Thus, if we denote the vector $v$ by $p_{\Lambda}(x)$ and
the vector $x - v$ by $q_{\Lambda}(x)$,
then we have (5).

(6)
Using (1), (2), (3) and (5), we can see
\begin{align*}
\Delta_{\Lambda}(x, x) & = \Delta_{\Lambda}(p_{\Lambda}(x), p_{\Lambda}(x))
+ 2 \Delta_{\Lambda}(p_{\Lambda}(x), q_{\Lambda}(x)) +
\Delta_{\Lambda}(q_{\Lambda}(x), q_{\Lambda}(x)) \\
& = \Delta_{\Lambda}(q_{\Lambda}(x), q_{\Lambda}(x)) =
\det (\Lambda) \langle q_{\Lambda}(x), q_{\Lambda}(x) \rangle.
\end{align*}
\QED

\begin{Corollary}
\label{cor:prop:elem:prop:B:lambda}
Let $f : V \to V'$ be a linear map of vector spaces over $\RR$, and let
$\langle\ ,\ \rangle$ and $\langle\ ,\ \rangle'$ be
inner products of $V$ and $V'$ respectively.
We assume that there is a positive constant $a$ with
$\langle f(x), f(x) \rangle' \leq a \langle x, x \rangle$ for
all $x \in V$.
Let $\Lambda = \{ v_1, \ldots, v_n \}$ be a set of linearly independent vectors in $V$, and
$\Lambda' = \{v'_1, \ldots, v'_{n'} \}$ a basis of $f(\Span(\Lambda))$.
Then, for all $x \in V$,
\[
\Delta_{\Lambda'}(f(x), f(x)) \leq a \frac{\det(\Lambda')}{\det(\Lambda)}\Delta_{\Lambda}(x, x).
\]
\end{Corollary}

\Proof
Let $x$ be an arbitrary element of $V$, and $x = v + y$ the decomposition of $x$
such that $v \in \Span(\Lambda)$ and $y$ is perpendicular to $\Span(\Lambda)$.
Then, by using (6) of Proposition~\ref{prop:elem:prop:B:lambda},
we can see that
\begin{align*}
a \frac{\det(\Lambda')}{\det(\Lambda)} \Delta_{\Lambda}(x, x) 
& = \det(\Lambda') a \langle y, y \rangle \\
& \geq \det(\Lambda') \langle f(y), f(y) \rangle' \\
& \geq \Delta_{\Lambda'}(f(y), f(y)).
\end{align*}
On the other hand, since $f(v) \in \Span(\Lambda')$, 
by (3) of Proposition~\ref{prop:elem:prop:B:lambda},
we can see
\[
\Delta_{\Lambda'}(f(x), f(x)) = \Delta_{\Lambda'}(f(y), f(y)).
\]
Thus, we get our corollary.
\QED

\setcounter{Theorem}{0}
\subsection{Proof}
\setcounter{Theorem}{0}
Let us begin with the following lemma.

\begin{Lemma}
\label{lemma:special:case:theorem:a}
Let $K$ be a finitely generated field over $\QQ$, and
$A$ an abelian variety over $K$.
Let $\Gamma$ be a subgroup of finite rank in $A(\overline{K})$.
Let $X$ be a subvariety of $A_{\overline{K}}$, and
$S$ an infinite subset of $X(\overline{K})$ with the following properties:
\begin{enumerate}
\renewcommand{\labelenumi}{(\arabic{enumi})}
\item
$S$ is generic, i.e., any infinite subsets of $S$ are Zariski dense in $X$.

\item
$S$ is small with respect to $\Gamma_{div} = \{ x \in A(\overline{K}) \mid
\text{$nx \in \Gamma$ for some positive integer $n$} \}$.
\end{enumerate}
Then, the stabilizer of $X$ in $A$ is positive dimensional.
\end{Lemma}

\Proof
First of all, since $S$ is infinite, $\dim(X) > 0$.
We fix a positive integer $n$ with $n > 2 \dim A$.
Enlarging $K$,
we may assume that $X$ is defined over $K$ and
there is a subgroup $\Gamma_0$ in $A(K)$
with $\Gamma_0 \subseteq \Gamma$ and $\Gamma_0 \otimes \QQ = 
\Gamma \otimes \QQ$.
By virtue of (1) of Proposition~\ref{prop:property:mini},
replacing $K$ by a finite extension field, $X$ by $[N](X)$, and
$S$ by an infinite subset of $[N](S)$, we may assume that
$(S, K)$ is $2$-minimized,
where $N$ is a positive integer.
Then, by virtue of Theorem~\ref{thm:equiv:cond:n:min},
there is an abelian subvariety $C$ of $A_{\overline{K}}$ such that
$\ZDD_2(S, K) = C$ and $\ZDD_n(S, K) = C^{n-1}$.

If $\dim C = 0$, then every element of $S$ is defined over $K$.
Here we use the following well known result, which is the special case of
Lang's conjecture:
\begin{quote}
``If $X(K)$ is Zariski dense in $X$,
then $X$ is a translation of an abelian subvariety of $A$.''
\end{quote}
Thus, $X$ is a translation of an abelian subvariety $G$ of $A$.
Then, $\Stab(X) = G$. Therefore, $\dim(\Stab(X)) = \dim G > 0$.

Next, we assume that $\dim(C) > 0$.
Let $\pi : A \to A/C$ be the natural homomorphism, and
$T = \pi(X)$. Let $X_T^n$ be the fiber product of $X$ over $T$ in $X^n$.
Then, we have a morphism $\beta_n : X_T^n \to A^{n-1}$.
Since $O_K(s)^n \subseteq X_T^n$, let
$Y$ be the Zariski closure of $\bigcup_{s \in S} O_K(s)^n$ in $X_T^n$.
Then, 
\[
\beta_n(Y) = \overline{\beta_n(Y)} \supseteq 
\overline{\beta_n\left( \bigcup_{s \in S} O_K(s)^n \right)}
= C^{n-1}.
\]
Therefore, we have
\[
\dim(X_T^n) \geq \dim(Y) \geq \dim (C^{n-1}).
\]
If the stabilizer of $X$ is finite, then
$\dim(X/T) \leq \dim(C) - 1$.
Thus,
\begin{align*}
\dim(X_T^n) - \dim (C^{n-1})  & =  (n\dim(X/T) + \dim(T)) - (n-1)\dim(C) \\
& \leq (n (\dim(C) - 1) + \dim(T)) - (n-1)\dim(C) \\
& = \dim(C) + \dim(T)  - n \\
& \leq 2 \dim(A) - n < 0.
\end{align*}
This is a contradiction. Therefore, $\dim(\Stab(X)) > 0$.
\QED

\bigskip
Let us start the proof of Theorem~\ref{thm:a:gen:lang:conj:fg:intro}.
We set $\Lambda = \{ \gamma_1, \ldots, \gamma_n \}$. Then, by using
the height pairing 
\[
\langle\ ,\ \rangle_L^{\overline{B}} : A(\overline{K}) \times
A(\overline{K}) \to \RR,
\]
we have the bilinear map
\[
\Delta_{\Lambda} : A(\overline{K})_{\RR} \times A(\overline{K})_{\RR}
\to \RR
\] 
as in \S\S\ref{subsec:pre:linear:algebra}.
Then, $\Delta_{\Lambda}(x,x) = \delta_L^{\overline{B}}(\gamma_1, \ldots, \gamma_n, x)$.

Let $\Stab(X)$ be the stabilizer of $X$ in $A$, and let
$\pi : A \to A' = A/\Stab(X)$ be the natural morphism.
We set $X' = \pi(X)$ and $\Gamma' = \pi(\Gamma)$.
Then, $\Stab(X')$ is trivial and $\pi^{-1}(X') = X$.
Let $L'$ be a symmetric ample line bundle on $A'$.
Then, by (2) of Proposition~\ref{prop:element:height:abelian}
(or Remark~\ref{rem:2:of:ineq:height}),
there is a positive number $a$ with
\[
\langle \pi(x), \pi(x) \rangle_{L'}^{\overline{B}} \leq a
\langle x, x \rangle_{L}^{\overline{B}}
\]
for all $x \in A(\overline{K})$.
Let $\Lambda' = \{ \gamma'_1, \ldots, \gamma'_{n'} \}$ be a basis
of $\Gamma' \otimes {\QQ}$.
Then, by Corollary~\ref{cor:prop:elem:prop:B:lambda},
\[
\Delta_{\Lambda'}(\pi(x), \pi(x)) \leq a \frac{\det(\Lambda')}{\det(\Lambda)}\Delta_{\Lambda}(x, x)
\]
for all $x \in A(\overline{K})$.
Thus, we can see that
the set $\{ x' \in X'(\overline{K}) \mid
\delta_{L'}^{\overline{B}}(\gamma'_1, \ldots, \gamma'_{n'}, x') \leq \epsilon \}$
is Zariski dense in $X'$ for every positive number $\epsilon$.
Here we assume that $\dim(X') > 0$. Then,
we can find a sequence $\{ x'_l \}_{l=1}^{\infty}$ in $X'(\overline{K})$
with the following properties:
\begin{enumerate}
\renewcommand{\labelenumi}{(\arabic{enumi})}
\item
If $l \not= m$, then $x'_l \not= x'_m$.

\item
$\{ x'_l \mid l=1,2, \ldots \}$ is generic in $X'$.

\item
$\delta_{L'}^{\overline{B}}(\gamma'_1, \ldots, \gamma'_{n'}, x'_l) < 1/l$
for all $l$.
\end{enumerate}
Here we claim the following.

\begin{Claim}
$\{ x'_l \mid l=1,2, \ldots \}$ is small with respect to $\Gamma'_{div}$.
\end{Claim}

In $A'(\overline{K}) \otimes \RR$, by (6) of Proposition~\ref{prop:elem:prop:B:lambda},

\[
\delta_{L'}^{\overline{B}}(\gamma'_1, \ldots, \gamma'_{n'}, x'_l) =
\Delta_{\Lambda'}(x'_l, x'_l) = \det(\Lambda') \langle x'_l - p_{\Lambda'}(x'_l),
x'_l - p_{\Lambda'}(x'_l) \rangle_{L'}^{\overline{B}} < 1/l.
\]
Here, since $\Gamma'_{\QQ}$ is dense in $\Gamma'_{\RR}$,
there is $y'_l \in \Gamma'_{\QQ}$ with
$\det(\Lambda') \langle x'_l - y'_l,
x'_l - y'_l \rangle_{L'}^{\overline{B}} < 1/l$.
Since $\Gamma'_{div}$ is a divisible group,
$y'_l$ comes from an element of $\Gamma'_{div}$, so that we may assume that
$y_l \in \Gamma'_{div}$.
Thus, if we set $z'_l = x'_l - y'_l$,
then $x'_l = y'_l + z'_l$, $y'_l \in \Gamma'_{div}$, and
$\det(\Lambda') \langle z'_l, z'_l \rangle_{L'}^{\overline{B}} < 1/l$.
Hence $\{ x'_l \mid l=1,2,\ldots \}$ is small with respect to $\Gamma'_{div}$
by (4) of Proposition~\ref{prop:element:prop:small}.

\medskip
By this claim together with Lemma~\ref{lemma:special:case:theorem:a},
we can see that $\dim (\Stab(X')) > 0$. This is a contradiction.
Therefore, $\dim(X') = 0$, say, $X' = \{ P' \}$.
Then, $\Delta_{\Lambda'}(P', P') \leq \epsilon$ for every $\epsilon > 0$.
Thus, $\Delta_{\Lambda'}(P', P') = 0$, which implies that $P' \in \Gamma'_{div}$.
Since $\pi : \Gamma_{div} \to \Gamma'_{div}$ is surjective,
there is $P \in \Gamma_{div}$ with $\pi(P) = P'$.
Then, $X = \Stab(X) + P$. Moreover, $\Stab(X)$ is an abelian
subvariety of $A$ because $X$ is a variety.
Thus, we get our theorem.

\begin{Remark}
Let $K$ be a finitely generated field over $\QQ$, $A$ an abelian variety over $K$,
and $X$ a geometrically irreducible subvariety of $A$.
Let 
\[
\langle\ , \ \rangle_L^{\overline{B}} :
A(\overline{K}) \times A(\overline{K}) \to \RR
\]
be the height pairing associated with a big polarization $\overline{B}$ and 
a symmetric ample line bundle $L$.
In the proof of this note,
we used only the following two fundamental results.

\medskip
$\bullet$ {\bf Bogomolov's conjecture over $\pmb{K}$}: 
If $\{ x \in X(\overline{K}) \mid \langle x, x \rangle_L^{\overline{B}} \leq \epsilon \}$
is Zariski dense in $X$ for every $\epsilon > 0$, then $X$ is a translation of
an abelian subvariety of $A$ by a torsion point.

\medskip
$\bullet$ {\bf Lang's conjecture over $\pmb{K}$ in the special case}:
If $X(K)$ is Zariski dense in $X$, then $X$ is a translation of an abelian subvariety of $A$.
\end{Remark}

\begin{Remark}
Even in the case where $K$ is a number field,
our proof is slightly simpler than Poonen's proof. For,
we avoid measure-theoretic arguments by considering a 
geometric trick.
\end{Remark}

\end{document}